\documentclass[letterpaper, 12 pt]{article}
\usepackage{epsfig, theorem, amsmath, amssymb, amsfonts}
\usepackage{doublespace}

\newcommand{\setleftmargin}[1]{
        \addtolength{\textwidth}{\oddsidemargin}
        \addtolength{\textwidth}{1in}
        \addtolength{\textwidth}{-#1}
        \setlength{\oddsidemargin}{-1in}
        \addtolength{\oddsidemargin}{#1}
        \setlength{\evensidemargin}{\oddsidemargin}
}
\newcommand{\setrightmargin}[1]{
        \setlength{\textwidth}{8.5in}
        \addtolength{\textwidth}{-\oddsidemargin}
        \addtolength{\textwidth}{-1in}
        \addtolength{\textwidth}{-#1}
}
\setleftmargin{1 in}
\setrightmargin{1 in}

\newcommand{\Z}{\mathbb{Z}}

\newcommand{\R}{\mathbb{R}}
\newcommand{\D}{\mathbb{D}}
\newcommand{\C}{\mathbb{C}}

\renewcommand{\l}{\xi}
\renewcommand{\Im}{\operatorname{Im}}

\newcommand{\ddt}{\frac{d}{dt}}
\newcommand{\atz}{\Big{|}_{t=0}}
\newcommand{\sumon}{\sum_{j=1}^n}

\DeclareMathOperator{\Sal}{Sal}
\newcommand{\subs}{\subseteq}
\newcommand{\YA}{\mathcal{Y}(\A)}
\newcommand{\ZA}{\mathcal{Z}(\A)}
\newcommand{\NA}{\mathcal{N}(\A)}
\newcommand{\MA}{\mathcal{M}(\A)}
\newcommand{\Rn}{\R^n}
\newcommand{\Rnd}{(\Rn)^{\vee}}
\newcommand{\Rs}{\R^\times}
\newcommand{\Cs}{\C^\times}
\newcommand{\Csn}{(\Cs)^n}
\newcommand{\Rts}{\R^2\smallsetminus\{0\}}
\newcommand{\Rd}{\R^{\vee}}
\newcommand{\Sala}{\Sal(\A)}
\newcommand{\Zt}{\Z_2}
\newcommand{\A}{\mathcal{A}}
\newcommand{\otn}{\{1,\ldots,n\}}

\newcommand{\rz}{\rho_z}

\newtheorem{theorem}{Theorem}[section]
\newtheorem{lemma}[theorem]{Lemma}

{\theorembodyfont{\rmfamily}
\theoremstyle{plain}

\newtheorem{example}[theorem]{Example}

\newtheorem{remark}[theorem]{Remark}
}

\newcommand{\qed}{\hfill \mbox{$\Box$}\medskip\newline}
\newenvironment{proof}{\noindent {\bf Proof:}}{\qed \par}

\newcommand{\<}{\left<}
\renewcommand{\>}{\right>}

\begin{document}
\begin{spacing}{1.1}

\noindent
{\LARGE \bf A nonhausdorff model for the complement of 
\smallskip\\ a complexified hyperplane arrangement}
\bigskip\\
{\bf Nicholas Proudfoot}\footnote{Partially supported
by an NSF Postdoctoral Research Fellowship}\\
Department of Mathematics, University of Texas,
Austin, TX 78712
\bigskip
{\small
\begin{quote}
\noindent {\em Abstract.}
Given a hyperplane arrangement $\A$ in a real vector space $V$ of dimension $d$,
we introduce a real algebraic prevariety $\ZA$, 
and exhibit the complement of the complexification
of $\A$ as the total space of an affine bundle over $\ZA$
with fibers modeled on the dual vector space $V^{\vee}$.
\end{quote}
}

\section{Statement}\label{results}
Let $V$ be a real vector space of dimension $d$, and let $f_1,\ldots,f_n$ be a collection
of nonconstant affine linear functions on $V$ such that the associated linear forms span the dual vector space $V^{\vee}$.
% For all $j\in\otn$, let $H_j=f_j^{-1}(0)$ and $H_j^\pm = f_j^{-1}(\R^\pm)$.
Let $\A$ denote the collection of affine hyperplanes $H_1,\ldots,H_n\subs V$,
where $H_j=f_j^{-1}(0)$ for all $j\in\otn$.
Let
$$\MA := V^\C\setminus\bigcup_{j=1}^n H_j^\C$$
be the complement of the complexification of $\A$.
Our goal is to introduce a new space which is which is naturally the base of a principal bundle
with structure group $V^{\vee}$ and total space $\MA$.
Note that a principal bundle with structure group equal to a vector space is the same
as an affine bundle whose associated vector bundle is trivial.

\begin{figure}[h]\label{npoints}
\centerline{\epsfig{figure=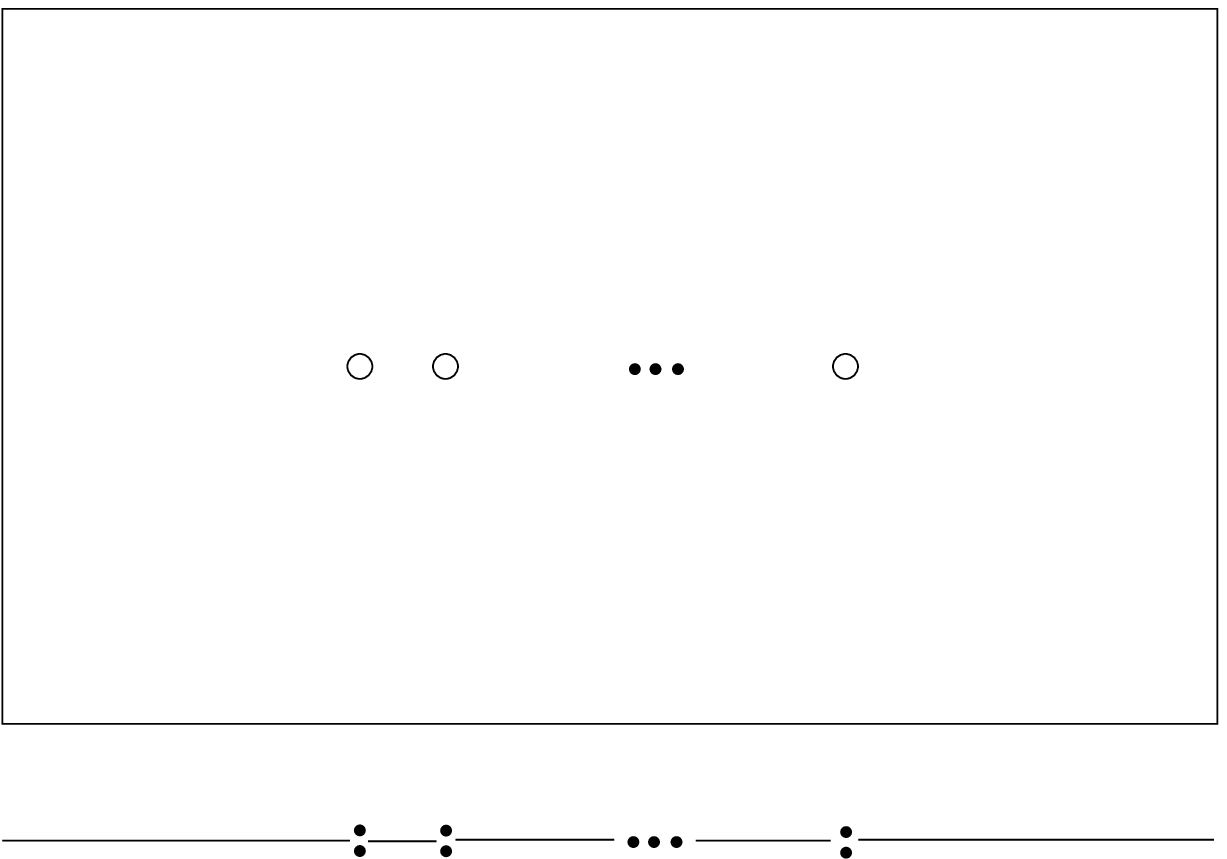, height=4cm}}
\caption{$\MA$ and $\ZA$, where $\A$ consists of $n$ distinct points on a line.}
\end{figure}

Our new space, which we call $\ZA$, will be a real algebraic prevariety of dimension $d$,
where ``prevariety'' means that $\ZA$ will not be Hausdorff in the analytic topology. 
Though this may sound nasty, our description will be quite simple.  
We construct $\ZA$ by gluing together
a collection of vector spaces $\{V_C\}$, each isomorphic to $V$, 
where the index $C$ ranges over the
{\em chambers} of $\A$, by which we mean the
connected components of the complement of $\A$ in $V$.  These vector spaces are attached to each
other along open sets, according to the following rule:
$$V_C\cap V_{C'} \cong V\,\,\setminus\!\!\bigcup_{\substack{\text{$H_j$ separates}\\ C \text{ from } C'}}\! H_j.$$
In Figure \ref{npoints} we illustrate the example of $n$ distinct points on a line;
here $\ZA$ is a real line with $n$ double points.
% For example, if $\A$ consists of $n$ distinct points on a line, then $\MA$ is isomorphic to $\Cs$,
% and $\ZA$ is a real line with a double point at the origin,
% as shown in Figure \ref{maza}.
More generally, $\ZA$ admits a natural map to $V$, and the number of points of $\ZA$ lying above a point
$q\in V$ is equal to the number of chambers $C$ which contain $q$ in their closures.

In this particular example, it is clear that there is a map from $\MA$ to $\ZA$ given by smooshing
the imaginary axis, and that this map is a weak homotopy equivalence (see Remark \ref{weak}).
In the following theorem, which is our main result, we strengthen this observation, and generalize it
to arbitrary arrangements.

\begin{theorem}\label{main}
The complexified complement $\MA$ may be realized as a principal bundle over $\ZA$ 
in the category of real analytic prevarieties,
with structure group naturally isomorphic to the additive group of the dual vector space $V^{\vee}$.
\end{theorem}

\begin{remark}\label{weak}
Since $\MA$ is isomorphic to the total space of a locally trivial fiber bundle over $\ZA$ with
contractible fibers, the projection from $\MA$ to $\ZA$ is a weak homotopy equivalence, i.e. it induces
isomorphisms on all homotopy and homology groups.  This map is not an honest homotopy equivalence
because it does not have a homotopy inverse; in particular, it admits no section.
\end{remark}

% \begin{remark}\label{vdual}
% In the analytic category, the multiplicative group $\Rpd$ is isomorphic via the logarithm map
% to the additive group $\Rd$, and a principal $\Rd$ bundle is a special case of an affine bundle of rank $d$.
% In our case, the structure group of the affine bundle 
% is naturally isomorphic to the vector space $V^{\vee}$.
% hence we may think of the group $\Rpd$ in Theorem \ref{main} as $\exp(V^*)$.
% \end{remark}

\begin{remark}
Theorem \ref{main} bears a strong similarity to the celebrated theorem of Salvetti \cite{Sa}
which exhibits a simplicial complex $\Sala$ that is homotopy equivalent to $\MA$.  
Salvetti's theorem offers some obvious advantages over Theorem \ref{main} -- the space in question
is a simplicial complex, one gets a homotopy equivalence rather than a weak homotopy equivalence,
and the Salvetti complex 
is an invariant of the oriented matroid associated to $\A$.
% can be associated to any pointed oriented matroid, regardless of whether
% it can be realized by a real hyperplane arrangements \cite{GR}.  
On the other hand, Theorem \ref{main}
has some advantages of its own.  First, the space $\ZA$ is quite easy to visualize, perhaps more so than $\Sala$.
Second, we are able to work real analyticlly, rather than simply topologically.
Finally, there is a canonical map from $\MA$ to $\ZA$,
while the homotopy equivalence in Salvetti's theorem requires some arbitrary choices.
\end{remark}

\begin{remark}\label{hypertoric}
In the special case where the arrangement $\A$ is central and
defined over the rational numbers, Theorem \ref{main}
may be interpreted in the world of hypertoric varieties, originally introduced by Bielawski and Dancer \cite{BD}.
A real hypertoric variety is a real variety of dimension $2d$ 
associated to a rational hyperplane arrangement
of rank $d$, and it carries an action
of the real algebraic gtorus $(\Rs)^d$.
% By taking the real locus we obtain a real variety of dimension $2d$ with an action of $(\Rs)^d$.
Let $\YA$ denote the open subset of this variety 
on which $(\Rs)^d$ acts freely.  One may show
that $\YA/\Zt^d$ is isomorphic to $\MA$, and $\YA/(\Rs)^d$ is isomorphic
to $\ZA$.  Hence $\MA$ is a principal bundle over $\ZA$ with structure group $(\Rs)^d/\Zt^d$, which is
analytically isomorphic to the additive group $\Rd$.
Though we will not use the language of hypertoric varieties in our proof of Theorem \ref{main},
those familiar with hypertoric varieties will observe that the proof is closely guided by this interpretation.
\end{remark}

\begin{remark}
Since $\ZA$ comes with an open cover $\{V_C\}$, it is natural to consider the Mayer-Vietoris
spectral sequence associated to this cover.  All multiple intersections of open sets are complements
of real hyperplane arrangements, and therefore have cohomology only in degree zero.  It follows that
the cohomology groups of $\ZA$ are equal to the homology groups of the complex $E_1^{0,*}$, where
$E_1^{0,q}$ is the direct sum of the $0$-th cohomology groups of all $(q+1)$-fold intersections.
Theorem \ref{main} tells us that the cohomology groups of $\ZA$ are isomorphic to those of $\MA$, which are
well understood by \cite{OS}.  It would be interesting to see the complex $E_1^{0,*}$ arise
in some independent combinatorial context.
\end{remark}

\begin{example}\label{smoosh}
We conclude the first section with a detailed analysis of the simplest case, that of a single point on a line,
which is the $n=1$ case in Figure \ref{npoints}.
Here $\MA$ is isomorphic to $\Cs$, and $\ZA$ is a real line with a double point at the origin.
We will write down explicitly the action of the additive group $\R$ on $\MA$, 
see that the orbit space is isomorphic to $\ZA$,
and check that the projection from $\MA$ to $\ZA$ is locally trivial.  
This will serve both as an illustration
of Theorem \ref{main}, and as an important tool to apply toward the proof 
of the general case in Section \ref{proof}.

First consider the action of $\Rs$ on $\Rts$ given by the formula 
\begin{equation}\label{action}
\lambda\cdot(x,y) = (\lambda x, \lambda^{-1} y).
\end{equation}
The map from $\Rts$ to $\R$ taking $(x,y)$ to $xy$ is surjective and $\Rs$-invariant.  The fiber of this map
over a nonzero number consists of a single $\Rs$ orbit, while the fiber over zero consists of two orbits, namely
the two coordinate axes.  
Hence the quotient of $\Rts$ by $\Rs$ is a line with a double point at the origin, which
we will call $\D$.  
The restriction of the action to the complement of either of the two coordinate axes in $\R^2$
induces a trivial $\Rs$ bundle over the quotient $\R$, 
therefore $\Rts$ is a principal $\Rs$ bundle over $\D$, trivialized
over the two copies of $\R$ in $\D$. 

Now consider the quotient of this entire picture by the subgroup $\Zt\subs\Rs$.  We now obtain an action
of $$\Rs\!/\Zt\cong\R_+\cong\R$$ on $$\big(\Rts\big)/\Zt\cong\Cs\!/\Zt\cong\Cs,$$ where the isomorphism
between $\R_+$ and $\R$ is given by the logarithm, and the isomorphism between $\Cs\!/\Zt$ and $\Cs$
is given by the map taking $\pm z$ to $iz^2$.
(We include the factor of $i$ for technical reasons to simplify notation in
the proof of the general case of Theorem \ref{main}.)  
Thus we obtain $\Cs$ as a principal $\R$ bundle over $\D$.
In fact, the structure group of the bundle is naturally dual to the original line,
hence we will denote the structure group $\Rd$.

For the sake of concreteness, let us write this action down in coordinates.
To make clear the distinction between the action of $\Rd$ on $\Cs$ and the action of $\Rs$ on $\Rts$ given
by Equation \eqref{action}, we will denote the $\Rd$ action with the symbol $*$.
Given an element $z \in \Cs$, choose a complex number $x+iy$ such that $z = i(x+iy)^2$.
Then for $\l\in\Rd$, 
% and $\pm\l\in\Rs/\Zt\cong\Rp$, 
we have
$$\l * z = i(e^{\l} x + e^{-\l} y i)^2 = -2xy + i(e^{2\l} x^2 - e^{-2\l} y^2) 
= \operatorname{Re}(z) + i(e^{2\l} x^2 - e^{-2\l} y^2).$$
Hence the action of $\Rd$ changes the imaginary part of $z$, and leaves the real part fixed.
If the real part of $z$ is nonzero, then $x$ and $y$ are both nonzero, and the orbit $\Rd * z$
is equal to the vertical real line through $z$.  On the other hand, if $z$ is purely imaginary, 
then either $x$ or $y$ is zero, and the orbit $\Rd * z$ is equal to the component
of the punctured imaginary axis containing $z$.
\end{example}

\begin{remark}
It is easy to become confused by the changes of variables in Example \ref{smoosh}.  We will
always adhere to the convention that when we use the letters $x$, $y$, and $z$, we have
$z = i(x+iy)^2$.  If we need letters to refer to the real and imaginary parts of $z$, as we will
in Lemma \ref{image}, we will use the notation $z = a+ib$.
\end{remark}

\section{Proof}\label{proof}
We now turn to the proof of the general case of Theorem \ref{main}.
Consider the affine linear map $$f:V\to \R^n$$ given by the functions $f_i$ which define the
hyperplanes of $\A$.
% , along with its linearization 
% $$\linf = f - f(0):V\to\Rn.$$
% % Let $$H_j^\pm = \{q\in V\mid \pm f_j(q)>0\}$$ and
% % $$K_j^\pm = \{q\in V\mid \pm \linf_j(q)>0\}.$$
The complexification $f^\C$ of $f$ 
induces a closed embedding of $\MA$ into $\Csn$.  
By Example \ref{smoosh},
$\Csn$ is a principal bundle over $\D^n$ with structure group $\Rnd$.
Let $\pi$
% :\Csn\to \D^n$ 
denote the projection from $\Csn$ to $\D^n$.

\begin{lemma}\label{image}
The image of $\MA$ under $\pi$ is isomorphic to $\ZA$.
\end{lemma}

\begin{proof}
Let us write $\D = \Rs \cup \{p,m\}$, where $p$ and $m$ stand for plus and minus.
Then the projection from $\Cs$ to $\D$ is given by sending $a+ib$ to $a$ if $a$ is nonzero,
and otherwise to $p$ or $m$, depending on the sign of $b$.

Let $\tilde a= (\tilde a_1,\ldots,\tilde a_n)$ be a point of $\D^n$ lying over a point $a=(a_1,\ldots,a_n)\in\Rn$.
Then $\tilde a$ lifts to an element $a+ib$ of $\MA$ if and only if 
$a$ lies in the image of $f$, and there exists $b=(b_1,\ldots,b_n)\in\Rn$
satisfying the following conditions:
\begin{eqnarray*}\label{conditions}
&(i)&  \text{$b+f(0)$ is in the image of $f$}\\
&(ii)&  \text{if $\tilde a_j = p$, then $b_j>0$}\\
&(iii)&  \text{if $\tilde a_j = m$, then $b_j<0$.}
\end{eqnarray*}

Consider the map $\tilde f_C:V_C\to \D^n$ lifting the map $f:V\to\Rn$, defined by the property
that if $f_j(q)=0$, then the $j^\text{th}$ coordinate of $\tilde f_C(q)$
is determined by the sign of $f_j$ on $C$.
These maps glue together to define an inclusion $$\tilde f:\ZA\to \D^n.$$
Furthermore, $\tilde a$ is in the image of $f_C$ if and only if $a$ is in the image of $f$ 
and $C$ is contained in the set
$$S(\tilde a) := \bigcap_{\tilde a_j = p}f_j^{-1}(\R_+)\,\,\cap\,\,\bigcap_{\tilde a_j = m}f_j^{-1}(\R_-).$$
It follows that $\tilde a$ is in the image of $\tilde f$ if and only if $a$ is in the image of $f$
and $S(\tilde a)$ is nonempty.
Suppose that $a = f(q_0)$ for some point $q_0\in V$.
Then the map $$q\mapsto f(q) - f(q_0) = f(q-q_0) - f(0)$$
takes $S(\tilde a)$ bijectively onto the set of points $b\in\Rn$ satisfying Conditions $(i), (ii),$ and $(iii)$
above, hence $S(\tilde a)$ is nonempty if and only if such a $b$ exists.
This proves that the image of $\tilde f$
coincides with the image of $\pi$ restricted to $\MA$.
\end{proof}

Let $W\subs\Rnd$ be the set of linear forms on $\Rn$ which are constant on the affine subspace $V$.
(Equivalently, $W$ is the kernel of the dual of the linearization of $f$.)
% We would like to think of $\Rpn$ as the connected algebraic group with Lie algebra $\Rnd$, hence
% we have an analytic isomorphism
% Consider the analytic isomorphism
% $$\exp:\Rnd\overset{\sim}\longrightarrow\Rpn$$ 
% given by identifying $\R$ with its dual in the natural way, and exponentiating each coordinate.
% Since $V^*$ is isomorphic to the quotient of $\Rnd$ by $W$, we will write
% $$\exp(V^*):=\Rpn/\exp(W).$$
% Thus $\exp(V^*)$ is the unique simply connected real algebraic group with Lie algebra $V^*$.
% It is abstractly isomorphic as an algebraic group to $\Rpd$,
% and analytically isomorphic to the additive group $V^*$ (see Remark \ref{vdual}). 
% 
Consider an element $z
% =(z_1,\ldots,z_n)
\in(\Cs)^n$, and choose $x,y\in\Rn$
such that $z_j = i(x+iy)^2$ for all $j\in\otn$. 
We define a function $\rz:W\to\R$ by the formula
$$\rz(\l) = |\!|e^\l\cdot(x,y)|\!|^2,$$
where $e^{\l}\in\exp(W)\subs(\Rs)^n$ acts coordinatewise
% on $\big(\Rts\big)^n$ 
as in Equation \eqref{action},
and the norm is 
% the restriction of 
the Euclidean norm on $\R^{2n}$.  Given a nonzero element $\xi\in W$, we compute
the directional derivative of $\rz$ along $\xi$:
\begin{equation}\label{deriv}
\ddt\big{|}\!\big{|}\exp(t\xi)\cdot(x,y)\big{|}\!\big{|}^2 =
\sumon\ddt\big{|}\!\big{|}(e^{t\xi_j}x_j,e^{-t\xi_j}y_j)\big{|}\!\big{|}^2 
= \sumon 2\,\xi_j\,\big(e^{2t\xi_j}x_j^2-e^{-2t\xi_j}y_j^2\big).
% &=& 2\<\xi,\Im(\exp(t\xi)*z)\>.
\end{equation}
Evaluating at $t=0$, we obtain 
\begin{equation}\label{pair}
\ddt\,\,\rz\big(\!\exp(t\xi)\big)\atz = \,\, 2\<\xi,\Im(z)\>,
\end{equation}
where $\<\cdot,\cdot\>$ is the natural pairing between $\Rnd$ and $\Rn$.

Let $$\NA = \pi^{-1}\big(\pi(\MA)\big) = \Rnd * \MA\subs\Csn.$$
By Lemma \ref{image}, $\NA$ is a principal bundle over $\ZA$ with structure group $\Rnd$.
We now prove a pair of lemmas that are analogous to the main theorem of \cite{KN}, which
lays the groundwork for the equivalence of quotients in symplectic and algebraic geometry.
This idea is central to the perspective alluded to in Remark \ref{hypertoric}.

\begin{lemma}\label{crit}
If $z\in\NA$, then $\l$ is a critical point of $\rz$ if and only if $e^\l * z\in\MA$.
\end{lemma}

\begin{proof}
Using the fact that $\rz(\l) = \rho_{\l * z}(\operatorname{id})$, we may immediately
reduce to the case where $\xi = 0$.
An element $z$ of $\NA$ automatically has its real part contained in the image of $f$,
hence it lies in $\MA$ if and only if $\Im(z)$ lies in the image of the linearization of $f$.
This in turn is the case if and only if the imaginary part of $z$ pairs trivially with every
$\xi\in W$ (this is how $W$ is defined).   
Equation \ref{pair} tells us that this happens if and only if zero
is a critical point of $\rz$.
\end{proof}

\vspace{-\baselineskip}
\begin{lemma}\label{unique}
If $z\in\NA$, then $\rz$ has a unique critical point.
\end{lemma}

\begin{proof}
By differentiating Equation \ref{deriv} we see that $\rz$
is convex, hence any critical point must be unique.  What remains is to prove existence.

Let us consider the behavior of the right hand side of Equation \ref{deriv} as $t$ becomes very large.
If $\xi_j>0$, then the $j^\text{th}$ term approaches positive infinity provided that $x_j\neq 0$,
otherwise it remains bounded.
Similarly, if $\xi_j<0$, then the $j^\text{th}$ term approaches positive infinity provided that $y_j\neq 0$,
and is otherwise bounded.  Hence the directional derivative of $\rz$ along $\xi$ at $t\xi*z$
is positive for large $t$ provided that there exists an index $j$ for which
either 
\begin{equation}\label{or}
\text{$\xi_j>0$\,\, and\,\, $x_j\neq 0$\hspace{1cm} or \hspace{1cm}$\xi_j<0$\,\, and\,\, $y_j\neq 0$.}
\end{equation} 
If this condition is satisfied, then when $t$ is large, the gradient of $\rz$ restricted to the
sphere of radius $t$ always points outward, hence 
$\rz$ must have a critical point somewhere on $\exp$ of the ball of radius $t$.
Hence it will suffice to prove that there exists an index $j$ satisfying Equation \ref{or}.

Recall that $z_j = i(x_j+iy_j)^2$.
If $x_j = 0$, then $z_j$ lies on the negative part of the imaginary axis;  if $y_j=0$, then it
lies on the positive part of the imaginary axis.
Since $z\in\NA$, there exists $w\in\MA$ lying in the same $\Rnd$ orbit as $z$. In particular, we have
$\Im(w_j) < 0$ whenever $x_j=0$, and $\Im(w_j)>0$ whenever $y_j=0$.
Suppose that Condition \eqref{or} fails for all $j$.  Then $\xi_j$ and $w_j$ have opposite signs whenever
$\xi_j\neq 0$, hence $\<\xi,w\> \neq 0$.  This contradicts the fact that $\xi\in W$ and $w\in \MA$.
\end{proof}

Lemmas \ref{crit} and \ref{unique} combine to tell us that
each $W$ orbit in $\NA$ contains a unique element of $\MA$.
Thus
$$\MA\cong\NA/W\cong\NA\times_{\Rnd} V^{\vee},$$  
which is the principal bundle over $\ZA$ induced from $\NA$ by the surjection
$$\Rnd\to\Rnd/W\cong V^{\vee}.$$
This completes the proof of Theorem \ref{main}.

\bigskip
\bigskip
\noindent {\em Acknowledgments.}
The author is grateful to Columbia University and the City of New York for their hospitality
during the writing of this paper.
\footnotesize{

}
\end{spacing}

\begin{thebibliography}{10}

\bibitem[BD]{BD}
R.~Bielawski and A.~Dancer.
\newblock The geometry and topology of toric hyperk\"ahler manifolds.
\newblock {\em Comm. Anal. Geom.} 8 (2000), 727--760.

% \bibitem[GR]{GR}
% I.M.~Gel$'$fand and G.L.~Rybnikov.
% \newblock Algebraic and topological invariants of oriented matroids.
% \newblock {\em Soviet Math. Dokl.} 40 (1990), no. 1, 148--152.

\bibitem[KN]{KN}
G.~Kempf and L.~Ness.
\newblock The length of vectors in representation spaces.  
\newblock {\em Algebraic geometry (Proc. Summer Meeting, Univ. Copenhagen, Copenhagen, 1978)},  pp. 233--243,
\newblock Lecture Notes in Math., 732, Springer, Berlin, 1979.

\bibitem[OS]{OS}
P.~Orlik and L.~Solomon.
\newblock Combinatorics and topology of complements of hyperplanes.
\newblock {\em Invent. Math.} 56 (1980), 167--189.

\bibitem[Sa]{Sa}
M.~Salvetti.
\newblock Topology of the complement of real hyperplanes in $\C^N$.
\newblock {\em Invent. Math.} 88 (1987), 603--618.

\end{thebibliography}
\end{document}